\documentclass{amsart}
\usepackage{amscd}
\usepackage{amssymb}

\def\k{k} 

\def\ker{\operatorname{ker}}
\def\nil{\operatorname{nil}}
\def\oo{\otimes}

\def\Spec{\operatorname{Spec}}

\def\tors{\operatorname{tors}}

\def\lra{\longrightarrow}
\newcommand{\mathdot}{{\mathbf{\scriptscriptstyle\bullet}}}
\newcommand{\Hcdh}{H_\cdh}
\def\red{\mathrm{red}}
\def\td{\mathrm{tr.\,deg}}
\newcommand{\comment}[1]{}  
\def \ie{{\it i.e.,\ }}

\def\cF{\mathcal F}

\def\HH{\mathcal{HH}}

\def\cO{\mathcal O}

\newcommand{\A}{\mathbb{A}}

\newcommand{\Q}{\mathbb{Q}}

\newcommand{\Z}{\mathbb{Z}}

\newcommand{\bbH}{\mathbb H}
\newcommand{\bbHcdh}{\bbH_{\cdh}}  

\def\cdh{\mathrm{cdh}}

\def\Cech{{\v C}ech\ }

\def\Carf{\operatorname{Carf}}

\def\nil{\operatorname{nil}}
\def\Pic{\operatorname{Pic}}

\def\map#1{{\buildrel #1 \over \lra}}

\def\onto{\twoheadrightarrow}
\def\defor{\overset\sim\onto}

\newcommand{\SchF}{\mathrm{Sch}/F}
\newcommand{\Schk}{\mathrm{Sch}/k}

\input xypic
\xyoption{all} 

\numberwithin{equation}{section}

\theoremstyle{plain}
\newtheorem{thm}[equation]{Theorem}
\newtheorem{cor}[equation]{Corollary}
\newtheorem{lem}[equation]{Lemma}
\newtheorem{prop}[equation]{Proposition}

\theoremstyle{definition}
\newtheorem{defn}[equation]{Definition}
\newtheorem{trick}[equation]{Standard Trick}
\theoremstyle{remark}
\newtheorem{rem}[equation]{Remark}
\newtheorem{ex}[equation]{Example}
\newtheorem{subremark}{Remark}[equation] 
\newtheorem{subex}[subremark]{Example} 

\begin{document}
\bibliographystyle{plain}

\title{Bass' $NK$ groups and $cdh$-fibrant Hochschild homology.}

\author{G. Corti\~nas}
\thanks{Corti\~nas' research was partially supported by Conicet 
and partially supported by grants PICT 2006-00836, UBACyT X051, 
PIP 112-200801-00900, and MTM2007-64704.}
\address{Dep. Matem\'atica, FCEyN-UBA\\ Ciudad Universitaria Pab 1\\
1428 Buenos Aires, Argentina}
\email{gcorti@gm.uba.ar}

\author{C. Haesemeyer}
\thanks{Haesemeyer's research was partially supported by NSF grant DMS-0652860}
\address{Dept.\ of Mathematics, University of California, Los Angeles CA
90095, USA}
\email{chh@math.ucla.edu}

\author{Mark E. Walker}
\thanks{Walker's research was partially supported by NSF grant DMS-0601666.}
\address{Dept.\ of Mathematics, University of Nebraska - Lincoln,
  Lincoln, NE 68588, USA}
\email{mwalker5@math.unl.edu}

\author{C. Weibel}
\thanks{Weibel's research was
supported by NSA grant MSPF-04G-184 and the Oswald Veblen Fund.}
\address{Dept.\ of Mathematics, Rutgers University, New Brunswick,
NJ 08901, USA} \email{weibel@math.rutgers.edu}

\date{\today}

\begin{abstract}
The $K$-theory of a polynomial ring $R[t]$ 
contains the $K$-theory of $R$ as a summand.
For $R$ commutative and containing $\Q$, we describe
$K_*(R[t])/K_*(R)$ in terms of Hochschild homology and the cohomology
of K\"ahler differentials for the $cdh$ topology.

We use this to address Bass' question, whether $K_n(R)=K_n(R[t])$
implies $K_n(R)=K_n(R[t_1,t_2])$. The answer to this
question is affirmative when $R$ is essentially of finite type over
the complex numbers, but negative in general.
\end{abstract}
\maketitle

\section*{}
\vskip-45pt

In 1972, H. Bass posed the following question
(see \cite{Bass73}, question (VI)$_n$):

\begin{center}
Does $K_n(R)=K_n(R[t])$ imply that $K_n(R)=K_n(R[t_1,t_2])$?
\end{center}

\noindent
One can rephrase the question in terms of Bass' groups $NK_n$,
introduced in \cite{Bass68}:
\begin{center}
Does $NK_n(R)=0$ imply that $N^2K_n(R)=0$?
\end{center}

More generally, for any functor $F$ from rings to an abelian category,
Bass defines $NF(R)$ as the kernel of the map $F(R[t])\to F(R)$ induced by
evaluation at $t=0$, and $N^2F=N(NF)$. Bass' question was inspired by
Traverso's theorem \cite{Traverso}, from which it follows that
$N\Pic(R)=0$ implies $N^2\Pic(R)=0$.

In this paper, we give a new interpretation of the groups 
$NK_n(R)$ in terms of Hochschild homology and the cohomology
of K\"ahler differentials for the $cdh$ topology, for commutative
$\Q$-algebras. This allows us to give a counterexample to Bass' question
in the companion paper \cite{TK} (see Theorem \ref{main00} below). 

To state our main structural theorem, recall from \cite{Weibel854} that
each $NK_n(R)$ has the structure of a module over the ring of 
big Witt vectors $W(R)$. It is convenient to use the countably 
infinite-dimensional $\Q$-vector spaces $t\Q[t]$ and $\Omega^1_{\Q[t]}$.
If $M$ is any $R$-module, then $M\oo t\Q[t]$ and $M\oo\Omega^1_{\Q[t]}$ 
are naturally $W(R)$-modules by \cite{DW}.

\begin{thm}\label{main01}
Let $R$ be a commutative ring containing $\Q$. Then there is a
$W(R)$-module isomorphism 
\[
N^2K_n(R)\ \cong \ \left( NK_n(R)\oo t\Q[t]\right) 
\ \oplus\ \left( NK_{n-1}(R)\oo\Omega^1_{\Q[t]}\right).
\]
Thus
$K_n(R)=K_n(R[t_1,t_2])$ iff $NK_n(R)=NK_{n-1}(R)=0$ iff $N^2K_n(R)=0$.

In addition, the following are equivalent for all $p>0$:
\begin{enumerate}
\item[a)] $K_n(R)=K_n(R[t_1,...,t_p])$
\item[b)] 
$
NK_n(R)=0 \textrm{ and }K_{n-1}(R)=K_{n-1}(R[t_1,...,t_{p-1}]).
$
\item[c)] $NK_q(R)=0$ for all $q$ such that $n-p<q\le n$.
\end{enumerate}
\end{thm}
\goodbreak


The equivalence of (a), (b) and (c) is immediate by induction,
using the formula for $N^2K_n$, 
and is included for its historical importance; see \cite{Vorst1}.
Theorem \ref{main01} also holds for the $K$-theory of schemes of 
finite type over a field; see \ref{N2Kdecomp} below.
\goodbreak

Theorem \ref{main01} allows us to reformulate Bass' question as follows:
\begin{equation*}
\text
{Does $NK_n(R) = 0$ imply that $NK_{n-1}(R) = 0$?}
\end{equation*}

\setlength{\leftmargini}{22pt}

\begin{thm}\label{main00}
a) For any field $F$ algebraic over $\Q$, the 2-dimensional normal algebra
$$R=F[x,y,z]/(z^2+y^3+x^{10}+x^7y)$$ 
has $K_0(R) = K_0(R[t])$ but $K_0(R) \neq K_0(R[t_1,t_2]).$

b) Suppose $R$ is essentially of finite type over a field 
of infinite transcendence degree over $\Q$. Then $NK_n(R) = 0$ implies
that $R$ is $K_n$-regular and, in particular, that $K_n(R) = K_n(R[t_1,t_2])$.
\end{thm}

Part (a) is proven in the companion paper \cite{TK}, using Theorem
\ref{main01}, while part (b) is proven below as Corollary \ref{cor:qBasslarge}.

The proof of Theorem \ref{main01} relies on methods developed in
\cite{chsw} and \cite{chw}, which allow us to compute the groups 
$NK_n$ and $N^pK_n$ in terms of 
the Hochschild homology of $R$, and of the $cdh$-cohomology of the
higher K\"ahler differentials $\Omega^p$, both relative to $\Q$.
The groups $NK_n(R)$ have a natural
bigraded structure when $\Q\subset R$, and it is convenient to take
advantage of this bigrading in stating our results. The bigrading
comes from the eigenspaces $NK_n^{(i)}(R)$ of the Adams operations $\psi^k$
(arising from the $\lambda$-filtration) and the eigenspaces of the 
homothety operations [r] (\ie base change for $t\mapsto rt$). 
This bigrading will be explained in Sections
\ref{sec:NHH} and \ref{sec:typical}; 
the general decomposition for Adams weight $i$ has the form:
\begin{equation}\label{bigrading}
 NK_n^{(i)}(R) \cong TK_n^{(i)}(R) \otimes_\Q t\Q[t].
\end{equation}
\noindent
Here $TK_n^{(i)}$ denotes the 
{\it typical piece} of $NK_n^{(i)}(R)$, defined as the
simultaneous eigenspace $\{ x\in NK_n^{(i)}(R): [r]x=rx,~ r\in R \}$.
(See Example \ref{ex:Carf-mod}.)
We provide a concrete description of the typical pieces in Theorem 
\ref{thm:NKformulas}, reproduced here:

\begin{thm}\label{main02}
If $R$ is a commutative $\Q$-algebra, then $NK_n^{(i)}(R)$ is
determined by its typical pieces $TK_n^{(i)}(R)$ and \eqref{bigrading}.
For $i\ne n,n+1$ we have:
\[ TK^{(i)}_{n}(R)\cong  \begin{cases}
HH^{(i-1)}_{n-1}(R)  &\text{if } i<n,\\
H_{\cdh}^{i-n-1}(R,\Omega^{i-1}) &\text{if } i\ge n+2.\end{cases}
\]
For $i=n,n+1$, we have an exact sequence:
\[0 \to TK_{n+1}^{(n+1)}(R) \to
\Omega^n_{R} \to H_{\cdh}^0(R,\Omega^n) \to TK_{n}^{(n+1)}(R) \to 0.
\]
\end{thm}

The special case $NK_0=\oplus NK_0^{(i)}$ of Theorem \ref{main02} is that 
for $R$ essentially of finite type over a field of characteristic zero,
with $d=\dim(R)$,

\begin{equation}\label{NK0}
NK_0(R) \cong \left( (R^+/R_{\red}) \oplus
\bigoplus\nolimits_{p=1}^{d-1}H_{\cdh}^p(R,\Omega^p)\right)
\otimes_\Q t\Q[t].
\end{equation}

Here $R^+$ is the seminormalization of $R_{\red}$;
we show in Proposition \ref{prop:seminorm} that $R^+=H^0_{\cdh}(R,\cO)$.
The dimension zero case of Theorem \ref{main02} is also revealing:

\begin{ex}\label{ex:nilp}
If $\dim(R)=0$ then we get
$NK_n(R)\cong HH_{n-1}(R,I)\otimes_\Q t\Q[t]$ for all $n$, where
$I$ is the nilradical of $R$. It is illuminating to compare this
with Goodwillie's Theorem \cite{Goodw}, which implies
that $NK_{n}(R)\cong NK_n(R,I)\cong NHC_{n-1}(R,I)$.
The identification comes from the standard observation
\eqref{NHCR} that the map $HH_*\to HC_*$ induces
$NHC_*(R,I)\cong HH_*(R,I)\otimes_\Q t\Q[t]$.
\end{ex}

The calculations of Theorem \ref{main02} 
for small $n$ are summarized in Table 1 when $\dim(R)=2$.
We will need the following cases of \ref{main02} in \cite{TK}, 
to prove Theorem \ref{main00}(a).
\goodbreak
%

\begin{thm}\label{thm:NS-intro}
Let $R$ be normal domain of dimension $2$ which is 
essentially of finite type over an algebraic extension of $\Q$. Then
\begin{enumerate}
\item[a)] $NK_0(R) = NK_0^{(2)}(R) \cong 
  H^1_{\cdh}(R,\Omega^1)\otimes_\Q t\Q[t]$ and
\item[b)] $NK_{-1}(R) = NK_{-1}^{(1)}(R) \cong 
  H^1_{\cdh}(R,\cO)\otimes_\Q t\Q[t].$
\end{enumerate}
\end{thm}

\begin{table} 
{\renewcommand{\arraystretch}{1.5}
\def\tabdots{\!\!$\cdots$\!\!} 
\begin{tabular}{ |c|c|c|c|c|c|c| }
\hline  & $i=1$ & $i=2$ & $i=3$ & $i=4$ & $i=5$ & $i=6$\hspace{-3pt} \\
\hline\hspace{-3pt} $TK_3^{(i)}(R)$ & 0 & $HH_2^{(1)}(R)$
& $\tors\Omega^2_R$\hspace{-3pt} & $\Omega^3_\cdh(R)/\Omega^3_R$\hspace{-2pt}
& $H_\cdh^1\Omega^4$\hspace{-2pt} & 0 \\
\hline\hspace{-3pt} $TK_2^{(i)}(R)$ & 0 & $\tors\Omega^1_R$
& $\Omega^2_\cdh(R)/\Omega^2_R$ & $H_\cdh^1\Omega^3$ & 0 \\
\cline{1-6}\hspace{-3pt} $TK_{1}^{(i)}(R)$ & $\nil(R)$
& $\Omega^1_\cdh(R)/\Omega^1_R$\hspace{-3pt} & $H_\cdh^1\Omega^2$ & 0 \\
\cline{1-5} $TK_0^{(i)}(R)$\hspace{-3pt}
& $R^+/R$ & $H_\cdh^1\Omega^1$ & 0\\
\cline{1-4}\hspace{-3pt} $TK_{-1}^{(i)}(R)$\hspace{-3pt}
& {\!}$H_\cdh^{1}\cO$\hspace{-4pt} & 0 \\
\cline{1-3} $TK_{-2}(R)$ & 0 \\
\cline{1-2}
\multicolumn{7}{c}{\parbox{4in}
{Table 1. The groups $TK_n^{(i)}(R)$ for $n\le3$, $\dim(R)=2$.}}
\end{tabular}
}\end{table}

Here is an overview of this paper: Section \ref{sec:NHH} reviews the
bigrading on the Hochschild and cyclic homology of $R[t]$
(and $X\times\A^1$), and Section \ref{sec:cdh} reviews the $cdh$-fibrant
analogue. Section \ref{sec:fibers} describes the sheaf cohomology of
the fibers $\cF_{HH}(X)$, $\cF_{HC}(X)$, etc.\ of $HH(X)\to\bbH_\cdh(X,HH)$,
etc. In Section \ref{sec:Bass} we use these fibers to prove
Theorem \ref{main01}, by relating $NK_{n+1}(X)$ to $H^{-n}\cF_{HH}(X)$.
We also show that Bass' question is negative for schemes
in Lemma \ref{lem:counterex}.

In Section \ref{sec:typical}, we give the detailed computations of the
typical pieces $TK_n^{(i)}(R)$ needed to 
establish \eqref{NK0} 
and Table 1;
these computations employ the main result of \cite{chwinf}.
In Section \ref{sec:qBasslarge}, we prove Theorem \ref{main00}(b),
that the answer to Bass' question is positive
provided we are working over a sufficiently large base field.
Finally, Section \ref{sec:surfaces} describes how Theorem \ref{thm:NS-intro}
changes if $R$ is of finite type over an arbitrary field of
characteristic~0: the map 
$NK_0(R)\to H^1_{\cdh}(R,\Omega^1_{/F})\oo_\Q t\Q[t]$
is onto, and an isomorphism if $NK_{-1}(R)=0$.

\vspace{-3pt}
\subsection*{Notation}
All rings considered in this paper should be assumed to be commutative
and noetherian, unless otherwise stated.
Throughout this paper, $k$ denotes a field of characteristic $0$ and
$F$ is a field containing $k$ as a subfield.
We write $\Schk$ for the category of separated 
schemes essentially of finite type over $k$.
If $\cF$ is a presheaf on $\Schk$, we write $\cF_\cdh$ for the associated
$cdh$ sheaf, and often simply write $H_\cdh^*(X,\cF)$ in place of
the more formal $H_\cdh^*(X,\cF_\cdh)$.

If $H$ is a functor on $\Schk$ and $R$ is an algebra essentially
of finite type, we occasionally write $H(R)$ for $H(\Spec R)$.
For example, $H^*_\cdh(R,\Omega^i)$ is used for $H^*_\cdh(\Spec R,\Omega^i).$
Note that, because the $cdh$ site is noetherian (every cover has a finite
subcovering) $H^*_\cdh$ sends inverse limits of schemes over 
diagrams with affine transition morphisms 
to direct limits. 
\goodbreak

If $H$ is a contravariant functor from $\Schk$ to spectra, 
(co)chain complexes, or abelian groups that takes
filtered inverse limits of schemes over 
diagrams with affine transition morphisms to colimits 
(as for example $K$, $HH$, $\bbH_\cdh(-,HH)$,
and $\cF_{HH}$), then for any $k$-algebra $R$, we abuse notation and
write $H(R)$ for the direct limit of the $H(R_\alpha)$ taken over all
subrings $R_\alpha$ of $R$ of finite type over $k$. (If $R$ is
essentially of finite type but not of finite type, the two definitions
of $H(R)$ agree up to canonical isomorphism.) 
In particular, we will use expressions like $\bbH_\cdh(R,HH)$ for
general commutative $\Q$-algebras even though we
do not define the $cdh$-topology for arbitrary $\Q$-schemes.

We use cohomological indexing for all chain complexes in this paper;
for a complex $C$, $C[p]^q = C^{p+q}$.
For example, the Hochschild, cyclic, periodic, and negative cyclic
homology of schemes over a field $k$
can be defined using the Zariski hypercohomology of
certain presheaves of complexes; see \cite{WeibelHC} and
\cite[2.7]{chsw} for precise definitions.  We shall write these
presheaves as $HH(/k)$, $HC(/k)$, $HP(/k)$ and $HN(/k)$, respectively,
omitting $k$ from the notation if it is clear from the context.

It is well known (see \cite[10.9.19]{WeibelHA94}) 
that there is an Eilenberg-Mac{\,}Lane functor
$C\mapsto|C|$ from chain complexes of abelian groups to spectra, and
from presheaves of chain complexes of abelian groups to presheaves of
spectra. 
This functor sends quasi-isomorphisms of complexes to weak
homotopy equivalences of spectra, and satisfies
$\pi_n(|C|)=H^{-n}(C)$.  For example, applying $\pi_n$ to the Chern
character $K\to|HN|$ yields maps $K_n(R)\to H^{-n}HN(R)=HN_n(R)$.  In
this spirit, we will use descent terminology for presheaves of
complexes.


\section{The bigrading on $NHH$ and $NHC$}\label{sec:NHH}

Recall that $k$ denotes a field of characteristic $0$. In this
section, we consider
the Hochschild and cyclic homology of polynomial extensions of
commutative $k$-algebras.  No great originality is claimed.
Throughout, we will use the chain level Hodge decompositions
$HH=\prod_{i\ge0} HH^{(i)}$ and $HC=\prod_{i\ge0} HC^{(i)}$.

The K\"unneth formula for Hochschild homology yields
\begin{equation}\label{NHHR}
NHH_n^{(i)}(R)\cong \left(HH_n^{(i)}(R)\otimes t\Q[t]\right) \oplus
\left(HH_{n-1}^{(i-1)}(R) \otimes\Omega^1_{\Q[t]}\right).
\end{equation}
 From the exact SBI sequence $0\to NHC_{n-1} \map{B} NHH_n \map{I}NHC_n\to0$
(see \cite[9.9.1]{WeibelHA94}), and induction on $n$, the map $I$ induces
canonical isomorphisms for each $i$:
\begin{equation}\label{NHCR}
NHC_n^{(i)}(R) \cong HH_n^{(i)}(R) \otimes t\Q[t].
\end{equation}

\begin{rem}\label{NHCX}
Both \eqref{NHHR} and \eqref{NHCR} generalize to non-affine
quasi-compact
schemes $X$ over $k$. Indeed, $NHH$ and $NHC$ satisfy Zariski descent
because $HH$ and $HC$ do and because, for any open cover $\{U_i\to X\}$, the
collection $\{U_i\times\A^1\to X\times\A^1\}$ is also a cover. Thus we have
\begin{align*}
NHH^{(i)}(X) & \cong  \bbH_{Zar}(X,NHH^{(i)}) \\
& \cong
\bbH_{Zar}(X,HH^{(i)})\otimes t\Q[t] \ \ \oplus\ \ 
\bbH_{Zar}(X,HH^{(i-1)})[1]\otimes\Omega^1_{\Q[t]} \\
& \cong  HH^{(i)}(X)\otimes t\Q[t] \ \
\oplus\ \ HH^{(i-1)}(X)[1]\otimes\Omega^1_{\Q[t]},
\end{align*}
and  $NHC^{(i)}(X)\!=\!\bbH_{Zar}(X,NHC^{(i)})\cong
\bbH_{Zar}(X,HH^{(i)})\otimes t\Q[t]\! \cong\! HH^{(i)}(X)\otimes t\Q[t]$.
\end{rem}

It is easy to iterate the construction $F\mapsto NF$. For example,
we see from \eqref{NHHR} and \eqref{NHCR} that
\begin{equation}\label{N2HCR}
N^2HC_n^{(i)}(R) \cong
\biggl(HH_n^{(i)}(R) \oo t\Q[t] \oo t\Q[t] \biggr) \oplus
\biggl(HH_{n-1}^{(i-1)}(R) \oo t\Q[t]\oo\Omega^1_{\Q[t]}\biggr).
\end{equation}
By induction, we see that 
$HH_{n-j}^{(i-j)}(R)\oo\bigl(t\Q[t]\bigr)^{\oo(p-j)} \oo
\bigl(\Omega^1_{\Q[t]}\bigr)^{\otimes j}$ will occur $\binom{p-1}j$ times
as a summand of $N^{p}HC_n^{(i)}(R)$ for all $j\ge0$. We may write
this as the formula:
\begin{equation}\label{NpHCR}
N^pHC_n^{(i)}(R) \cong \bigoplus_{j=0}^{p-1}
HH_{n-j}^{(i-j)}(R)\otimes_k \wedge^jk^{p-1}\oo\bigl(t\Q[t]\bigr)^{\oo(p-j)}
 \otimes \bigl(\Omega^1_{\Q[t]}\bigr)^{\otimes j}.
\end{equation}

\subsection*{Cartier operations on $NHH$ and $NHC$}

Let $W(R)$ denote the ring of big Witt vectors over $R$; it is well
known that in characteristic 0 we have $W(R)\cong\prod_1^\infty R$.
(See \cite[p.\,468]{Weibel854} for example.) 
Cartier showed in \cite{Cart} that the endomorphism ring $\text{Cart}(R)$
of the additive functor underlying $W$ consists of column-finite sums
$\sum V_m[r_{mn}]F_n$, using the {\it homotheties} $[r]$ (for $r\in R$),
and the Verschiebung and Frobenius operators $V_m$ and $F_m$.
Restricting the sum to $m\ge m_0$ yields a descending sequence of ideals
of $\text{Cart}(R)$, making it complete as a topological ring; 
$W(R)$ is the complete 
topological subring of all sums $\sum V_m[r_{m}]F_m$; see \cite{Cart}.

We will be interested in the intermediate (topological) subring 
$\Carf(R)$ of all row and column-finite sums $\sum V_m[r_{mn}]F_n$.
As observed in \cite[2.14]{DW}, there is an equivalence between the
category of $R$-modules and the category of continuous $\Carf(R)$-modules
given by the constructions in the following example.
(A left module $M$ is {\it continuous} if the annihilator ideal of each
element is an open left ideal.)

\begin{ex}\label{ex:Carf-mod}
If $M$ is any $R$-module, $N=M\otimes t\Q[t]$ is a 
continuous $\Carf(R)$-module (and hence a $W(R)$-module) via the formulas:
\[ [r]t^i=r^i t^i, \quad V_m(t^i)=t^{mi}, \quad F_m(t^i)=
\begin{cases} mt^{i/m} &\text{if } m|i,\\ 0 &\text{else.}\end{cases}
\]
The ring $W(R)=\prod_1^\infty R$ acts on $M\otimes t\Q[t]$ by
$(r_1,...,r_n,...)*\sum m_i t^i = \sum(r_im_i) t^i$.
Conversely, every continuous $\Carf(R)$-module $N$ has a ``typical
piece'' $M$, defined as 
the simultaneous eigenspace $\{ x\in N: [r]x=rx,~ r\in R \}$,
and $N\cong M\oo t\Q[t]$. 
\end{ex}

Recall that we can define operators $[r]$ on $NHH_n(R)$ and
$NHC_n(R)$, associated to the endomorphisms $t\mapsto rt$ of
$R[t]$. There are also operators
$V_m$ and $F_m$, defined via the ring inclusions $R[t^m]\subset R[t]$
and their transfers.  These operations commute with the Hodge
decomposition.
The following result
follows immediately from \cite[4.11]{DW} using
the observation that everything commutes
with Adams operations. 

\begin{prop}\label{prop:Carf-mod}
The operators $[r]$, $V_m$ and $F_m$ make each $NHC_n^{(i)}(R)$
into a continuous $\Carf(R)$-module, and hence a $W(R)$-module.
The $R$-module $HH_n^{(i)}(R)$ is its typical piece, and the canonical
isomorphism $NHC_n^{(i)}(R) \cong HH_n^{(i)}(R)\otimes t\Q[t]$ of
\eqref{NHCR} is an isomorphism of $\Carf(R)$-modules, the module structure
on the right being given in Example \ref{ex:Carf-mod}.
\end{prop}
\goodbreak

A similar structure theorem holds for $NHH_n(R)$ and its Hodge components,
using \eqref{NHHR}.
However, it uses a non-standard $R$-module structure on the typical piece
$HH_n(R)\oplus HH_{n-1}(R)$; see \cite[3.3]{DW} for details.

\smallskip
\begin{subremark}\label{rem:Carf-mod}
The conclusions of Proposition \ref{prop:Carf-mod} still hold for 
$NHC_n^{(i)}(X)$ and $HH_n^{(i)}(X)$ when $X$ is any scheme, where
$W(R)$ and $\Carf(R)$ refer to the ring $R=H^0(X,\cO)$.
That is, $HH_n^{(i)}(X)$ is an $R$-module and $NHC^{(i)}_n(X)$ is a
continuous $\Carf(R)$-module, isomorphic to $HH_n^{(i)}(X)\otimes t\Q[t]$.

This scheme version of \ref{prop:Carf-mod} is not stated
in \cite{DW}, which was written before the cyclic homology of schemes
was developed in \cite{WeibelHC}. 
However, the proof in \cite{DW} is easily adapted. Since
the operators $V_m$, $F_m$ and $[r]$ are defined on the underlying chain
complexes in \cite[4.1]{DW}, they extend to operations on the Hochschild and
cyclic homology of schemes. The identities required to obtain continuous
$\Carf(R)$-module structures all come from the K\"unneth formula for the
shuffle product on the chain complexes (see \cite[4.3]{DW}), so they also
hold for the homology of schemes.
\end{subremark}

\medskip
\section{$cdh$-fibrant $HH$ and $NHC$}\label{sec:cdh}

Now fix a field $F$ containing $k$; all schemes will lie in the category
$\SchF$ (essentially of finite type over $F$), in order to use
the $cdh$ topology on $\SchF$ of \cite{SVBK}.
All rings will be commutative $F$-algebras; because they are filtered direct
limits of finitely generated $F$-algebras, 
we can consider their $cdh$-cohomology.

If $C$ is any (pre-)sheaf of cochain complexes on $\SchF$, we can
form the $cdh$-fibrant replacement $X\mapsto\bbHcdh(X,C)$ and
write $\bbH_{\cdh}^n(X,C)$ for the $n$th cohomology of this
complex. (The fibrant replacement is taken with respect to the 
local injective model structure, as in \cite[3.3]{chsw}.)
For example, the $cdh$-fibrant replacement of a $cdh$ sheaf $C$
(concentrated in degree zero) is just an injective resolution, and
$\bbH_{\cdh}^n(X,C)$ is the usual cohomology of the $cdh$ sheaf
associated to $C$.

Hochschild and cyclic homology, as well as differential forms,
will be taken relative to $k$.
For $C=HH^{(i)}$, it was shown in \cite[Theorem 2.4]{chw} that
\begin{equation}\label{H-HH}
\bbHcdh(X,HH^{(i)})\cong\bbHcdh(X,\Omega^i)[i].
\end{equation}
This has the following consequence for $C=NHH^{(i)}$ and $NHC^{(i)}$.

\begin{lem} \label{lem:H-NHH}
Let $H^{(i)}$ denote either $HH^{(i)}$ or  $HC^{(i)}$, taken
relative to a subfield $k$ of $F$. Then
$\bbHcdh(X\times\A^1,H^{(i)})=\bbHcdh(X,H^{(i)}) \oplus \bbHcdh(X,NH^{(i)})$,
and:
\begin{align*}
\bbHcdh(X,NHH^{(i)}) & \cong \left(\bbHcdh(X,\Omega^i)[i] \otimes t\Q[t]
\right) \oplus
\left(\bbHcdh(X,\Omega^{i-1})[i]\otimes\Omega^1_{\Q[t]} \right); \\
\bbHcdh(X,NHC^{(i)}) & \cong \bbHcdh(X,\Omega^i)[i]\otimes t\Q[t].
\end{align*}
\end{lem}

\begin{proof}
The displayed formulas follow from \eqref{NHHR}, \eqref{NHCR} and
\eqref{H-HH}, using the first assertion and the fact that 
$- \otimes t\Q[t]$ commutes with $\bbHcdh$.
Thus it suffices to verify the first assertion.
By resolution of singularities, we may assume that $X$ is smooth.

Recall from \cite[3.2.2]{chsw} that the restriction of the $cdh$
topology to $\mathrm{Sm}/k$ is called the $scdh$-topology.
The product of any $scdh$ cover of $X$ with $\A^1$ is an
$scdh$ cover of $X\times\A^1$, and both $HH^{(i)}$ and
$HC^{(i)}$ satisfy $scdh$-descent by \cite[Thm. 2.4]{chw}.
Now by Thomason's Cartan-Leray Theorem \cite[1.56]{AKTEC} we have
\[ \bbHcdh(X\times\A^1,H^{(i)}) \cong \bbHcdh(X,H^{(i)}(-\times\A^1))
\cong\bbHcdh(X,H^{(i)}) \oplus \bbHcdh(X,NH^{(i)}).
\]
This gives the first assertion.
Alternatively, we may prove the first assertion by induction on $\dim(X)$,
using the definition of scdh descent to see that 
for smooth $X$ we have $H^{(i)}(X)=\bbHcdh(X,H^{(i)})$ and
\[
\bbHcdh(X\times\A^1,H^{(i)})=H^{(i)}(X\times\A^1)=H^{(i)}(X)\oplus NH^{(i)}(X).
\]
In particular, the first assertion holds when $\dim(X)=0$.
\end{proof}
\goodbreak

\begin{subremark}\label{rem:not-eft}
If $R$ is any commutative $F$-algebra,
the formulas of Lemma \ref{lem:H-NHH} hold for $X=\Spec(R)$ by naturality. This
is because we may write $R=\varinjlim R_\alpha$, where $R_\alpha$
ranges over subrings of finite type over $F$, and
$\bbHcdh(X,-)=\varinjlim\bbHcdh(\Spec(R_\alpha),-)$.
\end{subremark}

\begin{cor}\label{cor:Hcdh=0}
If $X=\Spec(R)$ is in $\SchF$, the modules
$\bbHcdh^n(X,HH^{(i)})$ and $\bbHcdh^n(X,NHC^{(i)})$ are zero
unless $0\le n+i < \dim(X)$ and $i\ge0$.

If $n\ge\dim(X)$ and $n>0$ then $\bbHcdh^n(X,HH)=0$.
\end{cor}

\begin{proof}
Because $\bbHcdh^n(X,\Omega^i)[i]=H^{i+n}_{\cdh}(X,\Omega^i)$,
this follows from \eqref{H-HH}, Lemma \ref{lem:H-NHH} and the fact
that $H^n_{\cdh}(X,\Omega^i)=0$ for $n\ge\dim(X)$, $n>0$. 
This bound is given in \cite[6.1]{chsw} for $i=0$, and in \cite[2.6]{chw}
for general $i$.
\end{proof}

Here is a useful bound on the cohomology groups appearing in
Lemma \ref{lem:H-NHH}. Given $X$, let $Q$ denote the total ring of fractions
of $X_\red$; it is a finite product of fields $Q_j$, and we let $e$ denote
the maximum of the transcendence degrees $\td(Q_j/k)$.

\begin{lem} \label{omegabound}
Let $X$ be in $\SchF$. If\ \ $i>e$ then $H^n_{\cdh}(X,\Omega^i)=0$ for all $n$.
\end{lem}

\begin{proof}
By \cite[12.24]{MVW}, we may assume $X$ reduced. 
Since we may write $X$ as an inverse limit of a sequence of affine
morphisms with the same ring of total factions $Q$, and $cdh$-cohomology
sends such an inverse limit to a direct limit, we may also assume that 
$X$ is of finite type over $F$. This implies that $e=\dim(X)+\td(F/k)$.

The result is clear if $\dim(X)=0$, since $H^n_\cdh(X, -) =
H^n_{Zar}(X, -)$ in that case.
Proceeding by induction on $\dim(X)$,
choose a resolution of singularities $X'\to X$ and observe that the
singular locus $Y$ and $Y\times_X X'$ have smaller dimension.
The hypothesis implies that $\Omega^i=0$ on $X'_{Zar}$, so
$H^n_{\cdh}(X',\Omega^i)=0$ by \cite[2.5]{chw}. The result now follows
by induction from the Mayer-Vietoris sequence of \cite[12.1]{SVBK}.
\end{proof}

If $R$ is a commutative ring, we write $R_{\red}$ and $R^+$ for the
associated reduced ring and the seminormalization of $R_{\red}$,
respectively. These constructions are natural with respect to localization,
so that we may form the seminormalization $X^+$ of $X_{\red}$ for any
scheme $X$. Because $X^+\to X$ is a universal homeomorphism,
we have $H^*_{\cdh}(X,-)\cong H^*_{\cdh}(X^+,-)$ for every $X$ in
$\Schk$, for any field $k$ of arbitrary characteristic.
The case $n=0$ with coefficients $\cO_\cdh$ is of special interest; 
recall our convention that $H^0_{\cdh}(X,\cO)$ denotes 
$H^0_{\cdh}(X,\cO_\cdh)$.

\begin{prop}\label{prop:seminorm}
For any algebra $R$, we have $H^0_{\cdh}(\Spec R,\cO) = R^+$.
Moreover,  for every $X$ in $\SchF$ we have
$H^0_{\cdh}(X,\cO) = \cO(X^+)$.
\end{prop}

\begin{proof}
We may assume $R$ and $X$ are reduced.
Writing $R=\varinjlim R_\alpha$ as in Remark \ref{rem:not-eft}, we have
$R^+=\varinjlim R_\alpha^+$ and
$H^0_\cdh(R,\cO)=\varinjlim H^0_\cdh(R_\alpha,\cO)$, so we may assume that
$R$ is of finite type. Thus the second assertion implies the first.
Since $H^0_\cdh(-, \cO)$ and $\cO(-^+)$ are Zariski sheaves, 
it suffices to consider the case when $X$ is affine.

Let $X = \Spec R$ be in $\SchF$, with $R$ reduced.
There is an injection $R\to Q$ with $Q$ regular
(for example, $Q$ could be the total quotient ring of $R$). By
\cite[6.3]{chsw}, $H^0_{\cdh}(\Spec Q,\cO)=Q$, so $R$ injects into
$H^0_{\cdh}(\Spec R, \cO)$.
This implies that $\cO_{\red}$ is a
separated presheaf for the $cdh$ topology on $\SchF$.  Thus, the ring
$H^0_{\cdh}(X, \cO)$ is the direct limit over all cdh-covers $p:U\to X$
of the \Cech $H^0$. (See \cite[3.2.3]{SGA4II}.)

Fix an element $b \in H^0_{\cdh}(\Spec R, \cO)$ and represent it by
$b \in \cO(U)$ for some $cdh$ cover $U \to X$.  Now recall from
\cite[12.28]{MVW} or \cite[5.9]{SVBK} that we may assume, by refining
the $cdh$ cover $U\to X$, that it factors as $U\to X'\to X$ where
$X'\to X$ is proper birational $cdh$ cover and $U\to X'$ is a Nisnevich cover.
If the images of $b \in \cO(U)$ agree in $U\times_X U$, \ie
$b$ is a \Cech cycle for $U/X$, 
then its images agree in $U \times_{X'} U$, \ie it is a
\Cech cycle for $U/X'$. But by faithfully flat descent, $b$ descends to
an element of $\cO(X')$.
Thus we can assume that $U$ is proper and birational over $X$.

Next, we can assume that the Nisnevich cover $p:U\to X$ is finite, surjective and birational.
Indeed, since $p$ is proper and birational we may consider the Stein
factorization $U \map{q} Y \map{r} X$. By
\cite[4.3]{EGAIII} or \cite[III.11.5 \& proof]{Hart}, 
$q_*(\cO_U)=\cO_Y$ and $r$ is finite surjective and birational. 
By  \cite[5.8]{SVBK},
$r$ is also a $cdh$ cover.  
Because $q_*(O_U)=O_Y$, the canonical
map $\cO_Y(Y)\to q_*(\cO_U)(Y)=\cO_U(U)$ is an isomorphism.
Hence $b$ descends to an element of $\cO(Y)$.
By Lemma \ref{lem:key},  $b$ lies in the seminormalization of $R$.
\end{proof}

\begin{lem}\label{lem:key} Let $A$ be a seminormal ring and $B$ a ring
between $A$ and its normalization.
Then the \Cech complex $A\to B \to B\otimes_A B$ is exact.
\end{lem}

\begin{proof}
We use Traverso's description of the seminormalization (see
\cite[p.\ 585]{Traverso}):
the seminormalization of a ring $A$
inside a ring $B$ is 
$$
A^+ = \{b\in B~\vert~ (\forall P\in \Spec A)~ b\in A_P+\mathrm{rad}(B_P)\}.
$$
Let $b\in B$ such that $1\otimes b = b\otimes 1$. We have to show that
$b\in A_P + \mathrm{rad}(B_P)$, for all primes $P$ of $A$. Let
$J = \mathrm{rad}(B_P)$; 
since $B_P/J$ is faithfully flat over the field $A_P/P$, the
image of $b$ in $B_P/J$ lies in $A_P/P$ by flat
descent. That is, $b\in A_P + J$, as required.
\end{proof}

\begin{rem} Even if $X$ is affine 
seminormal, it can happen that $H^i_{\cdh}(X,\cO) \neq 0$ for
some $i > 0$. For example, if $R$ denotes the subring $F[x,g,yg]$
of $F[x,y]$ for $g=x^3-y^2$ then it is easy to show that $R$ is seminormal
and that $H^1_{\cdh}(\Spec(R),\cO)=F$, because the normalization of $R$
is $F[x,y]$ and the conductor ideal is $gF[x,y]$.
For another example,
the normal ring of Theorem \ref{main00} has $H^1_{\cdh}(X,\cO)\ne0$,
by Theorems \ref{main01} and \ref{thm:NS-intro}(b). 
\end{rem}

\section{The fibers $\cF_{HH}$ and $\cF_{HC}$}\label{sec:fibers}

If $C$ is a presheaf of complexes on $\SchF$, we write $\cF_C$ for the
shifted mapping cone of $C\to\bbHcdh(-,C)$, so that
we have a distinguished triangle:
\begin{equation}\label{Ftriangle}
\bbHcdh(X,C)[-1] \to \cF_C(X)\to C(X)\to\bbHcdh(X,C)
\end{equation}

\begin{subex}\label{ex:O}
When $C$ is concentrated in degree $0$ we have $H^n\cF_C=0$ for all $n<0$.
For $C=\cO$ and $X=\Spec(R)$,
we see from Proposition \ref{prop:seminorm}
that $H^0\cF_{\cO}(X)=\nil(R)$,
$H^1\cF_{\cO}(X)=R^+/R$,
and $H^n\cF_{\cO}(X)=H_{\cdh}^{n-1}(X,\cO)$ for $n\ge2$.
Note that, if $X=\Spec R\in\SchF$, then $H^n\cF_{\cO}(X)=0$
for $n>\dim(X)$ by \cite[6.1]{chsw}.
\end{subex}

We now consider the Hochschild and cyclic homology complexes,
taken relative to a subfield $k$ of $F$.
For legibility, we write $\cF^{(i)}_{HH}$ for $\cF_{HH^{(i)}}$, etc.
By the usual homological yoga, $\cF_{HH}$ is the direct sum of the
$\cF^{(i)}_{HH}$, $i\ge0$, and similarly for $\cF_{HC}$.

\smallskip
\begin{subex}\label{ex:smoothFHH}
If $X$ is smooth over $F$ then $\cF_{HH} (X)\simeq 0$ by
\cite[2.4]{chw}.
\end{subex}

\smallskip
Lemma \ref{lem:H-NHH} and Remarks \ref{rem:not-eft} and \ref{NHCX}
imply the following analogue for $N\cF$.

\begin{lem}\label{lem:NFHCi}
If $X$ is in $\SchF$, or if $X=\Spec(R)$ for an $F$-algebra $R$,
we have quasi-isomorphisms:
\begin{align*}
N\cF_{HH}^{(i)}(X) & \cong \biggl(\cF_{HH}^{(i)}(X)\otimes t\Q[t]\biggr)\oplus
    \biggl(\cF_{HH}^{(i-1)}(X)[1]\otimes\Omega^1_{\Q[t]}\biggr);\\
N\cF_{HC}^{(i)}(X) & \cong \cF_{HH}^{(i)}(X)\otimes t\Q[t].
\end{align*}
\end{lem}

Mimicking the argument that establishes \eqref{N2HCR} and \eqref{NpHCR} yields:

\begin{cor}\label{cor:N2FHCi}
If $X$ is in $\SchF$, or if $X=\Spec(R)$ for an $F$-algebra $R$,
\[
N^2\cF_{HC}^{(i)}(X)\cong \bigl(\cF_{HH}^{(i)}(X)\oo t\Q[t]\oo t\Q[t]\bigr)
\oplus \bigl(\cF_{HH}^{(i-1)}(X)[1]\otimes t\Q[t]\oo\Omega^1_{\Q[t]}\bigr)
\]
and
\[
N^p\cF_{HC}^{(i)}(X)\cong \bigoplus_{j=0}^{p-1}
\cF_{HH}^{(i-j)}(X)[j]\otimes_k \wedge^jk^{p-1}\oo t\Q[t]^{\oo(p-j)}
\oo \bigl(\Omega^1_{\Q[t]}\bigr)^{\oo j}.
\]
\end{cor}

The cohomology of the typical pieces $\cF_{HH}^{(i)}(R)$
is given as follows.

\begin{lem}\label{lem:FHH} If $R$ is an $F$-algebra and $i\ge0$,
then there is an exact sequence:
\[ 0\to H^{-i}\cF_{HH}^{(i)}(R) \to \Omega_R^i \to H^0_{\cdh}(R,\Omega^i)
    \to H^{1-i}\cF_{HH}^{(i)}(R) \to 0. \]
For $n\ne i, i-1$ we have:
\[ H^{-n}\cF_{HH}^{(i)}(R) \cong \begin{cases}
HH_n^{(i)}(R) &\text{if } i<n,\\
H_{\cdh}^{i-n-1}(R,\Omega^i) &\text{if } i\ge n+2.\end{cases}
\]
\end{lem}

\begin{proof} As in Remark \ref{rem:not-eft}, we may assume $R$ is of finite type.
Since $HH_i^{(i)}(R)=\Omega^i_R$ for all $i\ge0$, and
$HH_n^{(i)}(R)=0$ when $i>n$ (see \cite[9.4.15]{WeibelHA94} or
\cite[4.5.10]{LodayHC92}), 
it suffices to use \eqref{H-HH} and to observe
that $\bbHcdh^{-n}(R,HH^{(i)})=H_{\cdh}^{i-n}(R,\Omega^i)$ vanishes
when $n>i$.
\end{proof}

\begin{ex}\label{ex:FHHR}
Let $X=\Spec(R)$ be in $\SchF$. Since $HH^{(0)}=\cO$, $\cF_{HH}^{(0)}(R)$
is described in Example \ref{ex:O}.
Applying Corollary \ref{cor:Hcdh=0} and Lemma \ref{lem:FHH} for $i>0$, and using
\cite[2.6]{chw} to bound the terms, we see that if $d=\dim(R)$ then
$H^n\cF_{HH}(X)=0$ for $n>d$. If $d=1$, then the only nonzero positive cohomology of
$\cF_{HH}$ is $H^1\cF_{HH}(R)=R^+/R$;
if $d>1$, we have:
\begin{align*}
H^1\cF_{HH}(R) \cong & \quad(R^+/R)\quad \oplus H_{\cdh}^1(X,\Omega^1)
\oplus \cdots \oplus H_{\cdh}^{d-1}(X,\Omega^{d-1}), \\
H^2\cF_{HH}(R) \cong &\ H_{\cdh}^1(X,\cO) \oplus
H_{\cdh}^2(X,\Omega^1)\oplus\cdots\oplus H_{\cdh}^{d-1}(X,\Omega^{d-2}),\\
\vdots\qquad & \qquad \vdots \\
H^d\cF_{HH}(R) \cong &\ H_{\cdh}^{d-1}(X,\cO).
\end{align*}
\end{ex}

\begin{ex}\label{ex:big-n}
When $R$ is essentially of finite type over $F$ and $\td(F/k)<\infty$,
$H^m\cF_{HH}(R)$ is Hochschild homology for large negative $m$.
To see this,
observe that $e=\td(R/k)$, the maximum transcendence degree of the
residue fields of $R$ at its minimal primes, 
is finite. 
Using Lemmas \ref{omegabound} and \ref{lem:FHH}, we get
$H^{-n}\cF_{HH}^{(i)}(R)=0$ and $H^{-n}\cF_{HH}^{(n)}(R)= \Omega^n_R$
for $i>n>e$, and hence
\[
H^{-n}\cF_{HH}(R) \cong HH_n(R) \text{ for all } n>e.
\]
\end{ex}

If $R=k\oplus R_1\oplus R_2\oplus\dots$ is graded,
and $\widetilde{HC}_*(R)=HC_*(R)/HC_*(k)$, it is well known that
the map $\widetilde{HC}_*(R)\map{S} \widetilde{HC}_{*-2}(R)$ is zero.
(See \cite[9.9.1]{WeibelHA94} for example.) 
In Lemma \ref{SBI-FHC} below, we prove a similar property for
$\cF_{HH}$ and $\cF_{HC}$, which we derive from Lemma \ref{lem:NFHCi}
using the following trick.

\begin{trick}\label{trick}
If $R$ is a non-negatively graded algebra, there is an algebra map
$\nu:R\to R[t]$ sending $r\in R_n$ to $rt^n$.
The composition of $\nu$ with evaluation at $t=0$
factors as $R\to R_0\to R$, and so
if $H$ is a functor on
algebras taking values in abelian groups, then the composition
$H(R)\map{\nu} H(R[t])\map{t=0} H(R)$
is zero on the kernel $\widetilde{H}(R)$ of $H(R)\to
H(R_0)$. 
Similarly,
the composition of $\nu$ with evaluation at $t=1$ is the identity. That is,
$\nu$ maps $\widetilde{H}(R)$ isomorphically onto a summand of $NH(R)$,
and $\widetilde{H}(R)$ is in the image of $(t=1):NH(R)\to H(R)$.
\end{trick}
\medskip

\begin{lem}\label{SBI-FHC}
If $R=\k\oplus R_1\oplus\cdots$ is a graded algebra, then for each $m$ the map
$\pi_m\cF_{HC}(R)\map{S}\pi_{m-2}\cF_{HC}(R)$ is zero and there
is a split short exact sequence:
\[
0 \to \pi_{m-1}\cF_{HC}(R) \map{B} \pi_m\cF_{HH}(R)
    \map{I} \pi_m\cF_{HC}(R) \to 0.
\]
Similarly, there are split short exact sequences:
\[
0\to \tilde \bbH_\cdh^{m+1}(R,HC) \map{B}
    \tilde \bbH_\cdh^m(R,HH)\map{I} \tilde \bbH_\cdh^m(R,HC) \to 0
\]
and
\begin{equation*}
0\to \tilde{\bbH}_\cdh^{m-1}(R,\Omega^{<i}) \map{B}
\tilde{H}_\cdh^{m-i}(R,\Omega^i) \map{I}
\tilde{\bbH}_\cdh^m(R,\Omega^{\le i}) \to 0.
\end{equation*}
\end{lem}

\begin{proof} It suffices to show that $I$ is onto and split.
By \cite[2.4]{chw}, $\cF_{HH}(\k)=\cF_{HC}(\k)=0$, so
$\tilde\cF_{HH}=\cF_{HH}$ and $\tilde\cF_{HC}=\cF_{HC}$.
By the standard trick \ref{trick}, 
it suffices to show that the maps
$N\pi_m\cF_{HH}(R)\to N\pi_m\cF_{HC}(R)$ and
$N\bbHcdh^m(R,HH)\to N\bbHcdh^m(R,HC)$ are split surjections.
But this is evident from the 
decompositions of 
$N\cF_{HC}^{(i)}(R)$ and $\bbHcdh(R,NHC^{(i)})$
in Lemmas \ref{lem:NFHCi} and \ref{lem:H-NHH}.

The third sequence is obtained from the second one by taking the
$i^{th}$ component in the Hodge decomposition,
described in Lemma \ref{lem:H-NHH}. 
\end{proof}
\begin{ex}\label{SBI-Homega}
\addtocounter{equation}{-1}
\begin{subequations}
Splicing the final sequences of Lemma \ref{SBI-FHC} together, we see that
the de Rham complexes are exact:
\begin{align} \label{dRH0}
0\to \k \to R & \map{d} \tilde{H}_\cdh^0(R,\Omega^1) \map{d}
\tilde{H}_\cdh^0(R,\Omega^2) \to\cdots\\ \label{dRHn}
0\to \Hcdh^n(R,\cO) & \map{d} \Hcdh^n(R,\Omega^1) \map{d}
    \Hcdh^n(R,\Omega^2) \to \cdots, \qquad n>0.
\end{align}
An analogous exact sequence
$$\cdots \to \pi_{m-1}\cF_{HH}(R)\map{d}
\pi_m\cF_{HH}(R)\map{d}\pi_{m+1}\cF_{HH}(R)\to \cdots$$ 
is obtained by splicing the other sequences in Lemma \ref{SBI-FHC}.
Using the interpretation of their Hodge components, described in
Lemma \ref{lem:FHH}, produces two more exact sequences:
\begin{align}\label{dRnil}
0\to \nil(R) \to & \tors\Omega^1_R \to \tors\Omega^2_R \to
\tors\Omega^3_R \to\cdots \\ \label{dROcdh}
0\to (R^+/R) \to & \Omega^1_\cdh(R)/\Omega^1_R \to
\Omega^2_\cdh(R)/\Omega^2_R \to \cdots.
\end{align}
Here we have written $\Omega^i_\cdh(R)$ for $\Hcdh^0(R,\Omega^i)$,
and $\tors\Omega^i_R$ is defined as the kernel of 
$\Omega^i_R\to\Omega^i_\cdh(R)$;
the notation reflects the fact that if $R$ is reduced then 
$\tors\Omega^i_R$ is the torsion submodule
of $\Omega^i_R$ (see Remark \ref{subrem:tors} below).
\end{subequations}
\end{ex}

\section{Bass' groups $NK_*(X)$}\label{sec:Bass}

In this section, we relate algebraic $K$-theory to our Hochschild
and cyclic homology calculations relative to the ground field $k=\Q$.
Consider the trace map
\[NK_{n+1}(X)\to NHC_n(X) = NHC_n(X/\Q)\]
induced by the Chern character. 
In the affine case, it is defined in \cite{WeibelNil};
for schemes it is defined using Zariski descent.
As explained in \cite{WeibelNil}, it arises from
the Chern character from the spectrum $NK(X)$ to the
Eilenberg-Mac{\,}Lane spectrum $|NHC(X)[1]|$ associated to the
cochain complex $NHC(X)[1]$. Note that our indexing conventions
are such that
$\pi_{n+1}|NHC(X)[1]|=H^{-n}NHC(X)=NHC_{n}(X).$

\begin{prop}\label{prop:NKdecomp}
Suppose that $R=\Gamma(X,\cO)$ for $X$ in $\SchF$, 
or that $X=\Spec(R)$ for an $F$-algebra $R$.
Then for all $n$, the Chern character induces a natural isomorphism
\begin{equation*}
NK_{n+1}(X)\cong H^{-n}\cF_{HH}(X) \otimes t\Q[t].
\end{equation*}
This is an isomorphism of graded $R$-modules, and even
$\Carf(R)$-modules, identifying the operations $[r]$, $V_m$ and
$F_m$ on $NK_*(X)$ with the operations on the right side described
in Example \ref{ex:Carf-mod}.
\end{prop}

\begin{proof} By Remark \ref{rem:not-eft}, we may suppose $X\in\SchF$.
By \cite[1.6]{chw}, the Chern character $K\to HN$
induces weak equivalences $\cF_K(X)\simeq |\cF_{HC}(X)[1]|$ and
$\cF_K(X\times\A^1)\simeq |\cF_{HC}(X\times\A^1)[1]|$.
Since for any presheaf of spectra $E$ we have a natural
objectwise equivalence $E(-\times\A^1) \simeq E \times NE$, we obtain
a natural weak equivalence
from $NK(X)$ to $|N\cF_{HC}(X)[1]|$. Now take homotopy groups and
apply Lemma \ref{lem:NFHCi}.

As observed in \cite[4.12]{DW}, the Chern character also commutes
with the ring maps used to define the operators $[r]$, $V_m$,
and with the transfer for $R[t^n]\to R[t]$ defining
$F_m$. That is, it is a homomorphism of $\Carf(R)$-modules. Since
the transfer is defined via the ring map $R[t]\to M_n(R[t^n])$,
followed by Morita invariance, there is no trouble in passing to
schemes.
\end{proof}

We now come to one of our main results, which implies 
Corollary \ref{main01}.

\begin{thm} \label{N2Kdecomp} For all $n$,
$N^2K_n(X) \cong \bigl(NK_n(X)\otimes t\Q[t] \bigr)
\oplus \bigl(NK_{n-1}(X)\otimes\Omega^1_{\Q[t]}\bigr)$, and
\[ N^{p+1}K_n(X) \cong \bigoplus_{j=0}^{p} NK_{n-j}(X)
\otimes \wedge^j\Q^{p}\otimes (t\Q[t])^{\otimes (p-j)} \otimes
\bigl(\Omega^1_{\Q[t]}\bigr)^{\otimes j}.
\]
This holds for every $X$ in $\SchF$, as well as for $\Spec(R)$ where
$R$ is an arbitrary commutative $F$-algebra.
\end{thm}

\begin{proof}
This is immediate from \ref{prop:NKdecomp} and Corollary \ref{cor:N2FHCi}.
\end{proof}

\begin{subremark}\label{rem:farrelljones}
Jim Davis has pointed out (see \cite{DavisNil})
that a computation equivalent to \ref{N2Kdecomp} can also be derived
--- for arbitrary rings $R$ ---
from the Farrell-Jones conjecture for the groups $\Z^r$. This particular
case is covered by F.\ Quinn's proof of hyperelementary assembly for
virtually abelian groups; see \cite{Quinn05}. 
\end{subremark}

As an immediate consequence of
\ref{N2Kdecomp} and \cite[XII(7.3)]{Bass68}, we deduce:
\goodbreak

\begin{cor}\label{abc} Suppose that $X$ is in $\SchF$,
or that $X=\Spec(R)$ for an $F$-algebra $R$. Then:
\begin{enumerate}
\item[a)] If $NK_n(X)=NK_{n-1}(X)=0$ then $N^2K_n(X)=0$.
\item[b)] If $NK_n(X)=0$ and $K_{n-1}(X)=K_{n-1}(X\times\A^{p})$ then
$K_n(X)=K_n(X\times\A^{p+1})$.
\item[c)] $K_n(X)=K_n(X\!\times\!\A^p)$ if and only if $NK_q(X)\!=0$
for all $q$ such that $n\!-\!p\!<\!q\!\le\! n$.
\end{enumerate}
\end{cor}
\goodbreak
Recall that $X$ is called {\it $K_n$-regular} if
$K_n(X)=K_n(X\times\A^p)$ for all $p$.

\begin{cor}\label{cor:Knregular}
Suppose that $X$ is in $\SchF$, 
or that $X=\Spec(R)$ for an $F$-algebra $R$.
Then the following conditions are equivalent:
\begin{enumerate}
\item[a)] $X$ is $K_n$-regular;
\item[b)] $NK_n(X)=0$ and $X$ is $K_{n-1}$-regular;
\item[c)] $NK_q(X)=0$ for all $q\le n$. 
\end{enumerate}
\end{cor}

\begin{subremark}
This gives another proof of Vorst's Theorem \cite[2.1]{Vorst1}
(in characteristic~0) that
$K_n$-regularity implies $K_{n-1}$-regularity, and extends it to
schemes.
\end{subremark}

The assumption that the scheme be affine is
essential in Bass' question ---
here is a non-affine example where the answer is negative.

\subsection*{Negative answer to Bass' question for nonaffine curves}

Let $X$ be a smooth projective elliptic curve over a number field
$k$ and let $L$ be a nontrivial degree zero line bundle with
$L^{\otimes3}$ trivial.  For example, if $X$ is the Fermat cubic
$x^3+y^3=z^3$, we may take the line bundle associated to the
divisor $P-Q$, where $P=(1:0:1)$ and $Q=(0:1:1)$.

\begin{lem}\label{lem:counterex}
Write $Y$ for the nonreduced scheme with the same underlying space
as $X$ but with structure sheaf $\cO_Y=\cO_X\oplus L = \mathrm{Sym}(L)/(L^2)$,
that is, $L$ is regarded as a square-zero ideal.

Then $NK_7(Y)=0$ but
$N^2K_7(Y)\cong NK_6(Y)\oo\Omega^1_{\Q[t]}$ is nonzero.
\end{lem}
\begin{proof}
In this setting, the relative Hochschild homology presheaf 
$HH_n(Y,L)$ is the kernel of $HH_n(Y)\to HH_n(X)$; sheafifying,
$\HH_n(Y,L)$ is the kernel of $\HH_n(Y)\to \HH_n(X)$.
Since $\Omega^1_X\cong\cO_X$ we see from Lemma~5.3 of \cite{chw}
that $\HH_n(Y,L)$ is: $L^{\otimes3}\oplus L^{\otimes5}$
if $n=4$; $L^{\otimes5}\oplus L^{\otimes5}$ if $n=5$; and
$L^{\otimes5}\oplus L^{\otimes7}$ if $n=6$. By Serre duality,
$H^*(X,L^{\otimes i})=0$ if $3\nmid i$ (cf.\ \cite[5.1]{chw}). By
Zariski descent, this implies that $HH_5(Y,L)\cong
H^1(X,\HH_4)\cong H^1(X,L^{\otimes3})\cong k$ and $HH_6(Y,L)=0$.
Since $\cF_{HH}(Y)\cong HH(Y,L)$, it follows from
\ref{prop:NKdecomp} and \ref{N2Kdecomp} that $NK_7(Y)=0$ but
$NK_6(Y)\cong t\Q[t]$ and 
$N^2K_7(Y)\cong NK_6(Y)\oo\Omega^1_{\Q[t]}\cong t\Q[t]\oo\Omega^1_{\Q[t]}.$
\end{proof}

We conclude this section by refining Proposition \ref{prop:NKdecomp}
and Corollary \ref{abc} to take account of the Adams/Hodge/%
$\lambda$-decompositions on K-theory and Hochschild homology, and
by establishing the triviality of $K_*^{(i)}(X)$ for $i\le0$.

Recall that by definition,
$K_n^{(i)}(X) = \{ x\in K_n(X)\otimes\Q:\psi^k(x)=k^ix\}$.
For $n < 0$, the Adams operations cannot be defined integrally.
However, it is possible to define the operations $\psi^k$ on
$K_n(X)\otimes\Q$ for $n<0$ using descending induction on $n$ and the
formula $\psi^k\{ x,t\}=k\{\psi^k(x),t\}$ in $K_{n+1}(X\times(\A^1-0))$
for $x\in K_n(X)$ and $\cO(\A^1-0)=F[t,1/t]$. 
This definition was pointed out in \cite[8.4]{Weibel-Pic}.

By \cite[2.3]{GW} or \cite[7.2]{chwinf},
the Chern character $NK_{n+1}(X)\to NHC_n(X)$
commutes with the Adams operations $\psi^k$ in the sense that
it sends $NK^{(i+1)}_{n+1}(X)$ to
$NHC_n^{(i)}(X)$ for all $i\le n$ (and to~0 if $i>n$).
Here is the $\lambda$-decomposition of
the isomorphism in Proposition \ref{prop:NKdecomp}:

\begin{prop}\label{prop:mainiso}
Suppose that $X\in\SchF$, or that $X=\Spec(R)$ for an $F$-algebra $R$.
Then for all $n$ and $i$, the Chern character induces a natural isomorphism:
\[  NK_{n}^{(i)}(X)   \cong  H^{1-n}\cF_{HH}^{(i-1)}(X) \otimes t\Q[t]. \]
In particular, if $i\le0$ then $NK_{n}^{(i)}(X)=0$ for all $n$.
\end{prop}
\goodbreak

\begin{proof}
By \cite{chwinf}, the Chern character $K\to HN$ sends
$K^{(i)}(X)$ to $HN^{(i)}(X)$.
The proof in \cite{chwinf} shows 
that the lift $\cF_K(X) \to \cF_{HN}(X)$, shown to be a
weak equivalence in \cite[1.6]{chw}, may be taken to send
$\cF_K^{(i)}(X)$ to $\cF_{HN}^{(i)}(X)$.
Since $HC\to HN$ sends $HC^{(i-1)}$ to $HN^{(i)}$\!,
the weak equivalence
$\cF_{HC}[1]\simeq\cF_{HN}$ identifies
$\cF_{HC}^{(i-1)}[1]$ and $\cF_{HN}^{(i)}$.
Finally $\cF_{HH}^{(i-1)}=0$ for $i\le0$.
\end{proof}

\begin{cor}
$K_n^{(i)}(X)\cong K_n^{(i)}(X\times\A^p)$ if and only if
$NK_{n-j}^{(i-j)}(X)=0$ for all $j=0,...,p-1$.
\end{cor}
\smallskip

\begin{thm}\label{thm:8.2}
For $X$ in $\SchF$ or $X=\Spec(R)$, and all integers $n$,
we have:
\begin{enumerate}
\item For $i<0$, $K_n^{(i)}(X)=0$.
\item For $i=0$, $K_n^{(0)}(X)\cong KH_n^{(0)}(X)\cong H^{-n}_\cdh(X,\Q)$.
\end{enumerate}
\end{thm}

\noindent Here $KH$ denotes the homotopy $K$-theory of \cite{WeibelKH}.
Theorem \ref{thm:8.2} answers Question 8.2 of \cite{Weibel-Pic}.

\begin{proof}
We first show that $K^{(i)}_n(X)\cong KH^{(i)}_n(X)$ when $i\le0$.
Covering $X$ with affine opens and using the Mayer-Vietoris sequences
of \cite[5.1]{WeibelKH}, it suffices to consider the case $X=\Spec(R)$.

Since $K(R)_{\Q}$ is the product of the eigen-components, the descent spectral
sequence $E^1_{p,q}=N^pK_q(R)_{\Q} \Rightarrow KH_{p+q}(R)_{\Q}$
(see \cite[1.3]{WeibelKH}) breaks up into one for each eigen-component.
If $i\le0$, the spectral sequence collapses by Proposition \ref{prop:mainiso}
to yield $K_n^{(i)}(R)\cong KH_n^{(i)}(R)$ for all $n$.

To determine the groups $KH_n^{(i)}(R)$ when $i\le0$, we use the $cdh$
descent spectral sequence of \cite[1.1]{HKH}.
If $i<0$, then the $cdh$ sheaf
$K^{(i)}_\cdh$ is trivial
as $X$ is locally smooth, so we have $KH^{(i)}_n(R)=0$ for all $n$.
If $i=0$ then the $cdh$ sheaf $K^{(0)}_\cdh$ is the sheaf $\Q_\cdh$;
see \cite[2.8]{sou}. Hence we
have $K^{(0)}_n(R)=KH^{(0)}_n(R)=H^{-n}_\cdh(X,\Q)$.
\end{proof}

\section{The typical pieces $TK_n^{(i)}(R)$}\label{sec:typical}

In this section, $R$ will be a commutative $F$-algebra.
The default ground field $k$
for K\"ahler differentials and Hochschild homology will be $\Q$.

As stated in \eqref{bigrading}, the Adams summands $NK_n^{(i)}(R)$ of
$NK_n(R)$ decompose as $NK^{(i)}_{n}(R)=TK^{(i)}_n(R)\otimes t\Q[t]$ for
each $n$ and $i$; the decomposition is obtained from an action of
finite Cartier operators precisely as the corresponding one for $NHC$
and $NHH$, explained in Section \ref{sec:NHH}.  The typical pieces
$TK^{(i)}_n(R)$ are described by the following formulas.

\begin{thm}\label{thm:NKformulas}
Let $R$ be a commutative $F$-algebra.
For $i\ne n,n+1$ we have:
\[ TK^{(i)}_{n}(R)\cong  \begin{cases}
HH^{(i-1)}_{n-1}(R),  &\text{if } i<n,\\
H_{\cdh}^{i-n-1}(R,\Omega^{i-1}) &\text{if } i\ge n+2.\end{cases}
\]

For $i=n,n+1$, the typical piece $TK^{(i)}_n(R)$
is given by the exact sequence:
\[ 0 \to TK_{n+1}^{(n+1)}(R) \to
\Omega^n_{R} \to H_{\cdh}^0(R,\Omega^n) \to TK_{n}^{(n+1)}(R) \to 0.
\]
\end{thm}
\goodbreak

\begin{proof} By Proposition \ref{prop:mainiso}, 
$TK_n^{(i)}=H^{1-n}\cF^{(i-1)}_{HH}$.
The rest is a restatement of Lemma  \ref{lem:FHH}.
\end{proof}

\begin{subremark}\label{subrem:trdeg}
If $R$ is essentially of finite type over a field $F$ whose
transcendence degree is finite over $\Q$, 
then the $TK_n^{(i)}(R)$ are finitely generated
$R$-modules. This fails if $\td(F/\Q)=\infty$ because then
$\Omega^i_{F/\Q}$ is infinite dimensional.
For instance, Example \ref{ex:nilp} implies that, for $R=F[x]/(x^2)$,
we have
$TK_2^{(2)}(R)=HH_1(R,x)=F\oplus\Omega^1_{F/\Q}.$
\end{subremark}
\goodbreak

\begin{cor}\label{NKneg}
Suppose that $R$ is essentially of finite type over $F$ and has
dimension $d$.
If $n<0$ then $NK_n^{(i)}(R)=0$ unless $1\le i\le d+n$,
in which case
\[ NK_n^{(i)}(R) = H_{\cdh}^{i-n-1}(R,\Omega^{i-1}) \oo t\Q[t]. \]
In particular, $NK_n(R)=0$ for all $n\le-d$.

If $d \ge2$ then:
\begin{align*}
NK_0(R) \cong & \quad\bigl[~(R^+/R)\hspace{5pt} \oplus H_{\cdh}^1(R,\Omega^1)
\oplus \cdots \oplus H_{\cdh}^{d-1}(R,\Omega^{d-1})\bigr]\oo t\Q[t], \\
NK_{-1}(R) \cong &\ \bigl[H_{\cdh}^1(R,\cO) \oplus
H_{\cdh}^2(R,\Omega^1)\oplus\cdots\oplus
H_{\cdh}^{d-1}(R,\Omega^{d-2})\bigr]\otimes t\Q[t],\\
\vdots\qquad & \qquad  \vdots \\
NK_{1-d}(R) \cong &\ H_{\cdh}^{d-1}(R,\cO) \otimes t\Q[t].
\end{align*}
If $d=1$ then $NK_0(R)=(R^+/R)\otimes t\Q[t]$ and $NK_n(R)=0$ for $n<0$.
\end{cor}

\begin{subremark}
The $d=1$ part of Theorem \ref{NKneg} holds for any
1-dimensional noetherian ring by \cite[2.8]{WeibelKA}.
\end{subremark}

\begin{subremark}\label{CHSW-bis}
Observe that \ref{NKneg} and \ref {cor:Knregular} imply that $R$ is
$K_{-d}$-regular. This
recovers the affine case of one of the main results in \cite{chsw}.
\end{subremark}

Here is a special case of the calculations in Theorem \ref{thm:NKformulas},
which proves Theorem \ref{thm:NS-intro}. We will use it to construct
the counterexample to Bass' question in 
the companion paper \cite{TK}.

\begin{thm}\label{thm:NS}
Let $F$ be a field of characteristic $0$ and $R$ a normal domain of
dimension $2$, essentially of finite type over $F$. Then
\begin{enumerate}
\item[a)] $H^1\cF_{HH}(R/F) \cong H^1_\cdh(R,\Omega^1_{/F})$,
\item[b)] $H^2\cF_{HH}(R/F) \cong H^1_\cdh(R,\cO)$,
\item[c)] $NK_0(R) \cong H^1_\cdh(R,\Omega^1)\oo t\Q[t]$, and
\item[d)] $NK_{-1}(R) \cong H^1_\cdh(R,\cO)\oo t\Q[t]$.
\end{enumerate}
\end{thm}

\begin{proof}
Parts (a) and (b) are immediate from Example \ref{ex:FHHR} 
and the
fact that $R$ is reduced and seminormal. Parts (c) and (d) follow from
(a) and (b) using Proposition \ref{prop:NKdecomp};
cf.\ Corollary \ref{NKneg}. 
\end{proof}
We introduce some notation to make the statement of the next theorem
more readable. The letter $e$ denotes the maximum transcendence degree of
the component fields in the total ring of fractions $Q$ of $R_\red$.
For simplicity, we write $\Omega_\cdh^i(X)$ for
$H^0_{\cdh}(X,\Omega^i)$, and we have written
$\Omega_{\cdh}^i(R)/\Omega^i_{R}$ for
the cokernel of $\Omega^i_{R}\to\Omega_{\cdh}^i(R)$.

\begin{defn}\label{def:EnR} 
For any commutative ring $R$ containing $\Q$, we define:
\begin{gather*} 
E_n(R)= \Omega_{\cdh}^n(R)/\Omega^n_R
\oplus\bigoplus\nolimits_{p=1}^{\infty} H^p_{\cdh}(R,\Omega^{n+p});
\\
\widetilde{HH}_n(R) = \ker\bigl( HH_{n}(R)\to\Omega^{n}_Q\bigr)
=\ker(\Omega^{n}_{R}\to\Omega^{n}_{Q})\oplus
\bigoplus\nolimits_{i=1}^{n-1}HH^{(i)}_{n}(R).
\end{gather*} 
\end{defn}

\begin{thm}\label{thm:05ab}
Let $R$ be a commutative ring containing $\Q$.
Then for all $n$:
\[
NK_{n}(R)\cong \bigl[ \widetilde{HH}_{n-1}(R)
\oplus E_n(R)\bigr]\otimes t\Q[t].
\]
\noindent
If furthermore $R$ is essentially of finite type over a field, and $n\ge e+2$,
then $NK_n(R)\cong HH_{n-1}(R)\otimes t\Q[t]$.
\end{thm}

\begin{proof}
Assembling the descriptions of the $TK_n^{(i)}(R)$ in Theorem
\ref{thm:NKformulas} yields the first assertion.
The second part is immediate from this and \ref{ex:big-n}.
\end{proof}
\goodbreak

\begin{subremark}
The Chern character
$NK_{n}(R)\to NHC_{n-1}(R)\cong HH_{n-1}(R)\oo t\Q[t]$
is an isomorphism for $n\ge e+2$. 
If $n\le e+1$, neither it nor the map
$H^{1-n}\cF_{HH}(R)\to HH_{n-1}(R)$ of  \ref{prop:NKdecomp} 
need be a surjection.
\end{subremark}
\smallskip

In order to compare the torsion submodules $\tors\Omega^*_R$ with the
typical pieces of $NK_*(R)$, we need the affine case 
of the following lemma. Following tradition, we write $F(X)$
for the total ring of fractions of $X_{\red}$.
That is, $F(X)$ is the product 
of the function fields of the irreducible components of $X_{\red}$.
When $X=\Spec(R)$ is affine, we write $Q$ instead of $F(X)$.

\begin{lem}\label{H0generic}
Let $X\in\SchF$; for $F(X)$ as above, the map
$\Omega_{\cdh}^i(X) \to \Omega^i_{F(X)}$ is an injection.
\end{lem}

\begin{proof} We may assume $X$ reduced, and proceed by induction on
$d=\dim(X)$, the case $d=0$ being trivial.
Choose a resolution of
singularities $X'\to X$ and let $Y$ be the singular locus of $X$,
with $Y'=Y\times_X X'$. By \cite[12.1]{SVBK}, there is a
Mayer-Vietoris exact sequence
\[
0\to \Omega_{\cdh}^i(X) \to \Omega_{\cdh}^i(X')\oplus\Omega_{\cdh}^i(Y)
\to \Omega_{\cdh}^i(Y') \map{\partial} H^1_{\cdh}(X,\Omega^i) \to \cdots.
\]
Since $F(Y)\subseteq F(Y')$, $\Omega^i_{F(Y)}\subseteq\Omega^i_{F(Y')}$.
Because $\dim(Y')<d$, the inductive hypothesis implies that
$\Omega^i_\cdh(Y)\to\Omega^i_\cdh(Y')$ is an injection.
Hence $\Omega^i_\cdh(X)\to\Omega^i_\cdh(X')$ is an injection.
But $X'$ is smooth, so by $scdh$ descent for $\Omega^i$
(see \cite[2.5]{chw}) we have
$\Omega^i_\cdh(X')\cong\Omega^i(X')\subset
\Omega^i_{F(X')} = \Omega^i_{F(X)}$.
\end{proof}

\begin{subremark}\label{subrem:tors}
Lemma \ref{H0generic} remains true if, instead of $\Omega^i$,
we use $\Omega^i_{/k}$ for $k\subseteq F$. In particular,
if $X=\Spec(R)$ is reduced affine, then
$\Omega^i_\cdh(R/k)=H^0_{\cdh}(R,\Omega^i_{/k})$ injects into
$\Omega^i_{Q/k}$. Thus $\tors(\Omega^i_{R/k})$, defined as the kernel of
$\Omega^i_{R/k}\to \Omega^i_\cdh(R/k)$ in \eqref{dRnil}, is
the torsion submodule of $\Omega^i_{R/k}$.
\end{subremark}

\begin{cor}\label{NKtors} For all $n\ge1$,
$TK_{n}^{(n)}(R)\cong\ker(\Omega^{n-1}_{R}\to \Omega^{n-1}_{Q})$.
\newline
In particular if $R$ is reduced, then 
$TK_{n}^{(n)}(R)$ is the torsion submodule of $\Omega^{n-1}_{R}$.
\end{cor}

\begin{proof} 
By \ref{thm:NKformulas}, $TK_{n}^{(n)}(R)$ is the kernel of
$\Omega^{n-1}_{R}\to \Omega^{n-1}_\cdh(R)$, so \ref{H0generic} applies. 
\end{proof}

The typical pieces of $NK^{(2)}_1(R)$ and $NK^{(2)}_2(R)$ of
\ref{thm:NKformulas} and \ref{NKtors} may be described as follows.

\begin{prop}\label{tau-tauF}
For all reduced $F$-algebras $R$, the typical pieces
$TK^{(2)}_1(R)=\Omega^1_{\cdh}(R)/\Omega^1_R$ and
$TK^{(2)}_2(R)=\tors(\Omega^1_R)$ fit into an exact sequence:
\begin{equation*}
0\to\tors(\Omega^1_R)\to\tors(\Omega^1_{R/F})\to\Omega^1_F\otimes(R^+/R)\to
\frac{\Omega^1_{\cdh}(R)}{\Omega^1_R} 
\to \frac{\Omega^1_{\cdh}(R/F)}{\Omega^1_{R/F}} \to 0.
\end{equation*}
\end{prop}

\begin{proof} We may assume $\Spec R\in\SchF$.
Recall from \cite[4.2]{chw} that there is a bounded second quadrant 
homological spectral sequence for all $p$ ($0\le i<p$, $j\ge0$):
\[
{}_p E^1_{-i,i+j}=\Omega^{i}_{F/k}\otimes_FHH^{(p-i)}_{p-i+j}(R/F)
\Rightarrow HH_{p+j}^{(p)}(R/k).
\]
When $p=1$, this spectral sequence degenerates 
to yield exactness of the bottom row in the following
commutative diagram; the top row is the First Fundamental Exact
Sequence for $\Omega^1$ \cite[9.2.6]{WeibelHA94}.
\[\xymatrix{ 
    &\Omega^1_F\otimes R\ar[d]\ar[r]&
    \Omega^1_R\ar[d]\ar[r]&\Omega^1_{R/F}\ar[d]\ar[r]&0\\
\qquad0\ar[r]&\Omega^1_F\otimes R^+\ar[r]&
\Omega^1_{\cdh}(R)\ar[r]& \Omega^1_{\cdh}(R/F)\ar[r]&0.
}\]
The upper left horizontal map is an injection because the left
vertical map is an injection. Now apply the snake lemma, using
Remark \ref{subrem:tors}.
\end{proof}

\section{Bass' question for algebras over large fields.}\label{sec:qBasslarge}

We will now show that the answer to Bass' question is positive for
algebras $R$ essentially of finite type over a field $F$ of
infinite transcendence degree over $\Q$.

Recall from Proposition \ref{prop:NKdecomp} that
$NK_{n+1}(R) \cong H^{-n}\cF_{HH}(R/\Q)\otimes t\Q[t].$
In light of this identification, the version 
of Bass' question stated before Theorem \ref{main00}
becomes the case $k=\Q$ of the following question:
\begin{equation}\label{q:HHBass}
\text
{Does $H^m\cF_{HH}(R/k)=0$ imply that $H^{m+1}\cF_{HH}(R/k)=0$?}
\end{equation}

In Theorem \ref{thm:qBassHHlarge}, we
show that the answer to question \eqref{q:HHBass}
is positive provided $R$ is of finite type over a field $F$ that has
infinite transcendence degree over $k$. The proof is essentially a
formal consequence of the K\"unneth formula in Lemma \ref{lem:bchange}.

\begin{lem}\label{lem:fibermodule}
Let $R$ be a commutative $F$-algebra,
and suppose $k$ is a subfield of $F$. 
Then $H^{-*}\cF_{HH}(R/k)$ and $\mathbb{H}^{-*}_\cdh(R,HH(/k))$
are graded modules over the graded ring $\Omega^{\mathdot}_{F/k}$.
\end{lem}

\begin{proof}
As in \ref{rem:not-eft}, we may suppose that $R$ is of finite type over $F$.
Consider the functor on $F$-algebras that associates to an $F$-algebra $A$
the Hochschild complex $HH(A/k)$. The shuffle product makes this into
a functor to $dg$-$HH(F/k)$-modules.
Since the $cdh$-site has a set of points (corresponding to valuations
by \cite[2.1]{GL}), 
we can use a Godement resolution to find a model for the
$cdh$-hypercohomology $\mathbb{H}_\cdh(-,HH(/k))$ which is also
a functor to $dg$-$HH(F/k)$-modules. It follows that there is a model for
$\cF_{HH}(R/k)$ that is a $dg$-$HH(F/k)$-module, functorially in $R$.
This implies the assertion, since 
$\Omega^\mathdot_{F/k}=H^{-\mathdot}HH(F/k)$.
\end{proof}

\goodbreak
\begin{lem}\label{lem:bchange} (K\"unneth Formula)
Suppose that $\Q\subseteq k \subseteq F_0 \subseteq F$ are fields.
Let $R_0$ be an $F_0$-algebra, and set
$R=F\otimes_{F_0}R_0$.
\begin{enumerate}
\item[i)] Let $T=\{t_i\}$ be transcendence basis of $F/F_0$;
writing $F[dT]$ for the exterior algebra on the set $\{dt_i\}$, we have
$\Omega^\mathdot_{F/F_0}=F[dT]$ and:
\[
\Omega^\mathdot_{F/k} \cong F[dT]\otimes_{F_0}\Omega^\mathdot_{F_0/k}
\]
In particular, the graded algebra homomorphism 
$\Omega^\mathdot_{F_0/k}\to \Omega^\mathdot_{F/k}$ is flat.

\medskip 

\item[ii)] $HH_*(R/k) \cong
\Omega^\mathdot_{F/k}\otimes_{\Omega^\mathdot_{F_0/k}}HH_*(R_0/k)
\cong F[dT]\otimes_{F_0}HH_*(R_0/k).$
\end{enumerate}
\end{lem}

\begin{proof} It is classical that $F[dT]=\Omega^\mathdot_{F/F_0}$.
The tensor product
decomposition of part i) follows from the fact that the fundamental sequence
\[
0\to F\otimes_{F_0}\Omega^1_{F_0}\to \Omega^1_{F}\to\Omega^1_{F/F_0}\to 0
\]
\goodbreak\noindent
is split exact. This proves i). To prove ii), choose a free chain
$dg$-$F_0$-algebra $\Lambda$ and a surjective quasi-isomorphism
of $dg$-algebras $\Lambda\defor R_0$. Then
$\Lambda'=F\otimes_{F_0}\Lambda\to F\otimes_{F_0}R_0=R$ is
a free chain model of $R$
as a $k$-algebra. Write $\Omega^\mathdot_{\Lambda}$ for differential forms;
consider $\Omega^\mathdot_{\Lambda}$ as a chain $dg$-algebra with the differential
$\delta$ induced by that of $\Lambda$. Note $\Lambda$
and $\Lambda'$ are homologically regular in the sense of \cite{cgg},
so that Theorem 2.6 of \cite{cgg} applies.
Combining this with part (i), we obtain
\begin{align*}
HH_*(R)=&HH_*(\Lambda')=H_*(\Omega^\mathdot_{\Lambda'})\\
       =&H_*(\Omega^\mathdot_F\otimes_{\Omega^\mathdot_{F_0}}\Omega^\mathdot_{\Lambda})=
    \Omega^\mathdot_F\otimes_{\Omega^\mathdot_{F_0}}H_*(\Omega^\mathdot_{\Lambda})\\
       =& ~\Omega^\mathdot_F\otimes_{\Omega^\mathdot_{F_0}}HH_*(R_0).\qedhere
\end{align*}
\end{proof}

Here is an easy consequence of Lemmas \ref{lem:fibermodule} and 
\ref{lem:bchange}.
\goodbreak

\begin{prop}\label{prop:fiber-bchange}
Suppose $\Q\subseteq k\subseteq F_0 \subseteq F$ are field extensions,
that $R_0$ is an $F_0$-algebra and $R = F\otimes_{F_0}R_0$.
Then there is an isomorphism of graded $\Omega^\mathdot_{F/k}$-modules
\[
F[dT]\oo_{F_0}
H^{-*}\cF_{HH}(R_0/k)
    \cong H^{-*}(\cF_{HH}(R/k)).
\]
\end{prop}

We also need the following lemma to prove the main result of this section.

\begin{lem}\label{lem:localvanish}
Let $R$ be essentially of finite type over $F\supset\Q$,
and let $H_n(R)$ denote either $HH_n(R)$ or $H^{-n}\cF_{HH}(R)$.
Assume that $H_{n_i}(R)=0$ for some finite set $\{n_1,\dots,n_r\}$
of positive integers.
Then there exist an $F$-algebra of finite type $R'$, and a 
multiplicatively closed set $S$ 
such that $R\cong S^{-1}R'$ and $H_{n_i}(R')=0$ for $1\le i \le r$.
\end{lem}

\begin{proof}
Because $R$ is essentially of finite type, it is the localization
$R=S^{-1}R''$ of some finite type $F$-algebra $R''$.  It is well known
that $HH_n(S^{-1}R'')\cong S^{-1}HH_n(R'')$ (see \cite[9.1.8]{WeibelHA94}),
and $H^{-n}\cF_{HH}(S^{-1}R'')\cong S^{-1}H^{-n}\cF_{HH}(R'')$ by
\cite[2.8--9]{chw}. 

Because $R''$ is of finite type over $F$, we may write 
$R''=F\otimes_{F_0}R_0$ for some finitely generated field extension 
$F_0$ of $\Q$ and some finite type $F_0$-algebra $R_0$. 
Note $R_0$ is essentially of finite type over $\Q$, whence
$H_p(R_0)$ is a finitely generated $R_0$-module ($p\ge 0$).
By \ref{lem:bchange} and/or \ref{prop:fiber-bchange},
$H_p(R'')$ is isomorphic, as an $R''$-module, to a
direct sum of copies of $R''\otimes_{R_0}H_q(R_0)$ 
with $q\le p$.  In particular, $M=\bigoplus_{i=1}^rH_{n_i}(R'')$ 
is a finite sum of $R''$-modules, each of which is a 
--- possibly infinite --- direct sum of copies of one finitely 
generated module. 

Given that $M$ has this form, the hypothesis that $S^{-1}M=0$ 
implies that there exists a nonzero
element $s\in Ann(M)\cap S$. Consider the finite type $F$-algebra 
$R'=R''[1/s]$.  Then $R\cong S^{-1}R'$ and we have
$\bigoplus_i H_{-n_i}(R')=M[1/s]=0$.
\end{proof}

\begin{thm}\label{thm:qBassHHlarge}
Suppose $k\subset F$ is an extension with $\td(F/k)=\infty$, and $R$
is essentially of finite type over $F$.  If $H^n(\cF_{HH}(R/k)) = 0$,
then $H^{m}(\cF_{HH}(R/k)) = 0$ for all $m\geq n$.
\end{thm}

\begin{proof}
By Lemma \ref{lem:localvanish}, we may assume that $R$ is of finite
type over $F$. There is a finitely generated field extension
$F_0\subset F$ of $k$ and a finite type $F_0$-algebra $R_0$ such that
$R = R_0\otimes_{F_0}F$. Note that $\td(F/F_0)=\infty$.
By Lemma \ref{lem:bchange} and Proposition \ref{prop:fiber-bchange},
$\Omega^i_{F/F_0}\otimes_{F_0} H^{n+i}(\cF_{HH}(R_0/k))$ 
is a direct summand of $H^n(\cF_{HH}(R/k))$ for each $i\geq 0$. 
Since $\Omega^i_{F/F_0}\ne0$ 
for all $i$, all the $H^{n+i}(\cF_{HH}(R_0/k))$ vanish as well.
Similarly, $H^{m}(\cF_{HH}(R/k))$ is a direct sum of copies of the groups
$\Omega^j_{F/F_0}\otimes_{F_0}H^{m+j}(\cF_{HH}(R_0/k))$ for $j\geq 0$,
all of which vanish when $m\ge n$, as we just observed.
\end{proof}
\smallskip

\begin{cor}\label{cor:qBasslarge}
Let $\Q\subset F$ be a field extension of infinite transcendence degree,
and suppose $R$ is essentially of finite type over $F$. Then $NK_n(R) = 0$
implies that $R$ is $K_n$-regular.
\end{cor}

\begin{proof}
Combine Theorem \ref{thm:qBassHHlarge} with 
Proposition \ref{prop:NKdecomp} and Corollary \ref{cor:Knregular}.
\end{proof}

Here is another proof of Corollary \ref{cor:qBasslarge}, 
which is essentially due to Murthy and Pedrini
and given in their 1972 paper \cite{MP}; they stated the result only for
$n\le1$ because transfer maps for higher $K$-theory and the
$W(R)$-module structure had not yet been discovered. We are grateful
to Joseph Gubeladze \cite{Gubeladze} for pointing this out to the authors.

\begin{lem}\label{MP1.4}
If $R$ is an algebra over a field $k$ of characteristic~0,
$N^pK_n(R[t])\to N^pK_n(R\otimes_k k(t))$ is injective.
\end{lem}

\begin{proof}
The proof in \cite[1.3--1.6]{MP} goes through, taking into account 
that the norm map and localization sequences used there for 
$K_0$, $K_1$ are now known for all $K_n$.
\end{proof}

\begin{lem}\label{MP2.1}
Suppose that $k$ is an algebraically closed field of 
infinite transcendence degree over $\Q$, 
and that $R$ is a finitely generated $k$-algebra.
If $NK_n(R)$ is zero, then 
$K_n(R)\map{\simeq} K_n(R[x_1,...,x_p])$ for all $p>0$.
\end{lem}

\begin{proof}
Muthy and Pedrini prove this in \cite[2.1.]{MP}; 
although their result is only stated for $i\le1$, 
their proof works in general. Note that since $NK_n(R)$ has the form
$TK_n(R)\oo t\Q[t]$ by \eqref{bigrading} (a result which was not
known in 1972), $NK_n(R)$ is torsionfree,
and has finite rank if and only if it is zero.
\end{proof}

\begin{proof}[Proof of Corollary \ref{cor:qBasslarge}]
Let $\Phi$ denote the functor $N^pK_n$. If $k\subset k_1$ is a finite
algebraic field extension and $R$ is a $k$-algebra, then
$\Phi(R)\to \Phi(R\oo_k k_1)$ is an injection because its composition
with the transfer $\Phi(R\oo_k k_1)\to \Phi(R)$ is multiplication by
$[k_1:k]$, and $\Phi(R)$ is a torsionfree group. Since $\Phi$ commutes
with filtered colimits of rings, $\Phi(R)\to \Phi(R\oo_k\bar{k})$
is an injection. Thus Lemma \ref{MP2.1} suffices to prove 
Corollary \ref{cor:qBasslarge} when $R$ is of finite type.
\end{proof}

\goodbreak
\section{$NK_0$ of surfaces}\label{sec:surfaces}

We conclude with a general description for affine surfaces of 
the canonical map $\Omega^1_{F}\oo_F NK_{-1}\to NK_0$. 
This sheds light on the difference between the cases of small and
large base fields, and also explains some results of
\cite{WeibelNorm}.

If $R$ is a 2-dimensional noetherian ring then $NK_0(R)$ is the direct
sum of $NK_0^{(1)}(R)=N\Pic(R)$ and $NK_0^{(2)}(R)$
\goodbreak

\begin{thm}\label{thm:Chuckoffer}
Let $R$ be a 2-dimensional normal domain of finite type over a field
$F$ of characteristic~0.
There is an exact sequence: 
\begin{multline*}
0\to NK^{(2)}_1(R)\to
\left({H^0(R,\Omega^1_{/F})}/{\Omega^1_{R/F}}\right)\oo t\Q[t]
\\  \to \Omega^1_{F}\otimes_F NK_{-1}(R)\to
NK_0(R)\to H^1_\cdh(R,\Omega^1_{/F})\oo t\Q[t]\to 0.
\end{multline*}
\end{thm}
\begin{proof}
Consider the following short exact sequence of sheaves in $(\SchF)_\cdh$:
\[
  0\to \Omega^1_F\otimes_F\cO\to \Omega^1\to \Omega^1_{/F}\to 0
\]
Applying $H_\cdh$ yields
\begin{multline*}
0\to \Omega^1_F\otimes_F R\overset\iota\to H^0(R,\Omega^1)\to
    H^0(R,\Omega^1_{/F})\overset\partial\to\\
\Omega^1_{F}\otimes_F H^1_\cdh(R,\cO)\to
    H^1_\cdh(R,\Omega^1)\to H^1_\cdh(R,\Omega^1_{/F})\to 0
\end{multline*}
Note that, because $\Omega^1_R\to \Omega^1_{R/F}$ is onto, the map
$\partial$ kills the image of $\Omega^1_{R/F}$. Similarly, the image
of $\iota$ is contained in that of $\Omega^1_R$. Thus we obtain
\begin{multline*}
0\to {H^0(R,\Omega^1)}/{\Omega^1_R}\to
    {H^0(R,\Omega^1_{/F})}/{\Omega^1_{R/F}}\to
\\  \Omega^1_{F}\otimes_F H^1_\cdh(R,\cO)\to
H^1_\cdh(R,\Omega^1)\to H^1_\cdh(R,\Omega^1_{/F})\to 0
\end{multline*}
Now apply $\oo t\Q[t]$ and use \ref{thm:NKformulas} and parts c) and d)
of \ref{thm:NS}.
\end{proof}
\begin{cor}
Let $R$ be a 2-dimensional normal domain of finite type over a field
$F$ of characteristic~0.
If $NK_{-1}(R)=0$ then $NK_0(R) \cong H^1_\cdh(R,\Omega^1_{/F})\oo t\Q[t]$.
\end{cor}

\begin{ex}\label{ex:bchangedim2}
Let $R$ be a 2-dimensional normal domain of finite type over $\Q$, and
put $R_F=R\oo F$. By \ref{prop:NKdecomp} and \ref{prop:fiber-bchange},
\begin{equation}\label{bchange}
NK_*(R_F)\cong NK_*(R)\oo\Omega^*_{F/\Q}.
\end{equation}
Keeping track of the $\lambda$-decomposition, as in \ref{thm:NKformulas},
we see from Theorem \ref{thm:NS-intro} that
\[
TK_1^{(2)}(R_F)\cong TK_1^{(2)}(R)\oo F \cong
H^0(R,\Omega^1)\!\oo\!F/\Omega^1_{R}\!\oo\!F \cong
H^0(R_F,\Omega^1_{/F})/\Omega^1_{R_F/F}.
\]
 From Theorem \ref{thm:Chuckoffer} we get an exact sequence
\begin{equation}\label{short}
0\to \Omega^1_{F/\Q}\otimes_F NK_{-1}(R_F)\to NK_0(R_F)\to
    H^1_\cdh(R_F,\Omega^1_{/F})\oo t\Q[t]\to 0
\end{equation}
Using \eqref{bchange} and \ref{thm:NS-intro} again,
we see that the sequence \eqref{short} is isomorphic to the sum
\begin{gather*}
(0\to \Omega^1_{F/\Q}\oo H^1_\cdh(R,\cO)\oo t\Q[t] \map{\simeq}
\Omega^1_{F/\Q}\oo H^1_\cdh(R,\cO)\oo t\Q[t]\to 0\to 0)\\\oplus\\
(0\to 0 \to F \otimes H^1_\cdh(R,\Omega^1)\oo t\Q[t]   \map{\simeq}
F\otimes H^1_\cdh(R,\Omega^1)\otimes t\Q[t]\to 0)
\end{gather*}
For example, for
$R_F:=F[x,y,z]/(z^2+y^3+x^{10}+x^7y)$
the results of \cite{TK} show that:
\begin{align*}
NK_{-1}(R_F)&=F\otimes t\Q[t]\\
 NK_0(R_F)&=\Omega^1_{F/\Q}\oo t\Q[t]\cong\bigoplus_{p=1}^{\td(F)} F\oo t\Q[t].
\end{align*}
In other words, both typical pieces $TK_{-1}(R_F)$ and $TK_0(R_F)$ are
$F$-vectorspaces, but while $\dim_FTK_{-1}(R_F)=1$ for all $F$, any
cardinal number $\kappa$ can be realized as $\dim_F TK_0(R_F)$ for an
appropriate $F$.
\end{ex}

\subsection*{Acknowledgements}
The authors would like to thank M.~Schlichting, whose contributions
go beyond the collaboration \cite{chsw}. We would also like to thank
W.~Vasconcelos, L.~Avramov, E.~Sell and J.~Wahl
for useful discussions.

\end{document}